\documentclass[12pt]{article}%
\usepackage{amssymb}%
\usepackage{amsmath}%
\usepackage{amsfonts}%
\usepackage{graphicx}
\providecommand{\U}[1]{\protect\rule{.1in}{.1in}}
\newtheorem{theorem}{Theorem}

\newtheorem{corollary}[theorem]{Corollary}

\newtheorem{lemma}[theorem]{Lemma}

\newtheorem{proposition}[theorem]{Proposition}

\begin{document}

\title{\textbf{Triangle Decompositions of Planar Graphs}}
\author{C. M. (Kieka) Mynhardt\thanks{Supported by a Discovery Grant from the Natural
Sciences and Engineering Research Council of Canada.}\ \ and Christopher M.
van Bommel\thanks{Supported by a Julie Payette Research Scholarship and an
Andr\'{e} Hamer Postgraduate Prize from the Natural Sciences and Engineering
Research Council of Canada.}\\Department of Mathematics and Statistics\\University of Victoria, P.O. Box 1700 STN CSC\\Victoria, BC, \textsc{Canada} V8W 2Y2\\{\small kieka@uvic.ca; cvanbomm@uvic.ca}}
\date{January 12, 2015}
\maketitle

\begin{abstract}
A multigraph $G$ is triangle decomposable if its edge set can be partitioned
into subsets, each of which induces a triangle of $G$, and rationally triangle
decomposable if its triangles can be assigned rational weights such that for
each edge $e$ of $G$, the sum of the weights of the triangles that contain $e$
equals $1$.

We present a necessary and sufficient condition for a planar multigraph to be
triangle decomposable. We also show that if a simple planar graph is
rationally triangle decomposable, then it has such a decomposition using only
weights $0,1$ and $\frac{1}{2}$. This result provides a characterization of
rationally triangle decomposable simple planar graphs. Finally, if $G$ is a
multigraph with $K_{4}$ as underlying graph, we give necessary and sufficient
conditions on the multiplicities of its edges for $G $ to be triangle and
rationally triangle decomposable.

\end{abstract}

\noindent\textbf{Keywords:\hspace{0.1in}}Planar graphs; Triangle
decompositions; Rational triangle decompositions

\noindent\textbf{AMS Subject Classification 2010:\hspace{0.1in}}05C10; 05C70

\section{Introduction}

We consider multigraphs, in which multiple edges between vertices are allowed,
but loops are not, and reserve the term \emph{graph} for a simple graph. For a
graph $H$, a multigraph $G$ is $H$-\emph{decomposable} if its edge set can be
partitioned into subsets, each of which induces a subgraph isomorphic to $H$.
Such a partition is called an $H$-\emph{decomposition} of $G$. A $K_{3}%
$-decomposition is also called a \emph{triangle decomposition}, and a $K_{3}%
$-decomposable multigraph is also said to be \emph{triangle decomposable}.
Given a multigraph $G$, a\emph{\ rational }$K_{3}$\emph{-decomposition} of $G$
is an assignment of nonnegative rational numbers, called \emph{weights}, to
the copies of $K_{3}$ in $G$ such that for each edge $e$ of $G$, the sum of
the weights of the triangles that contain $e$ equals $1$. If $G$ admits a
rational $K_{3}$-decomposition, we say that $G$ is \emph{rationally triangle
decomposable }or \emph{rationally }$K_{3}$-\emph{decomposable}.

We present a necessary and sufficient condition for a planar multigraph to be
triangle decomposable. This result implies that a maximal planar graph is
$K_{3}$-decomposable if and only if it is Eulerian. We also present results on
rationally $K_{3}$-decomposable planar multigraphs, including a
characterization of rationally $K_{3}$-decomposable planar simple graphs.

Triangle decompositions of graphs have a long history, beginning with the
following problem raised by W. S. B. Woolhouse in 1844 in \emph{The Lady's and
Gentleman's Diary} \cite[as cited by Biggs in \cite{Biggs}]{Wool}:

\begin{center}
\textquotedblleft Determine the number of combinations that can be made of $n
$ symbols,\\[0pt]$p$ symbols in each; with this limitation, that no
combination of $q$ symbols\\[0pt]which may appear in any one of them shall be
repeated in any other.\textquotedblright
\end{center}

\noindent A version of this problem (in which each pair of symbols appears
\emph{exactly} once) was solved for $p=3$ and $q=2$ by Thomas Kirkman \cite[as
cited in \cite{Biggs}]{Kirk} in 1847. Structures satisfying these constraints
became known as Steiner triple systems in honour of Jakob Steiner \cite[as
cited in \cite{Biggs}]{Stein}, who independently posed the question of their existence.

Simple necessary conditions for a connected multigraph $G$ to be triangle
decomposable are that $G$ be Eulerian and $|E(G)|\equiv0\ (\operatorname{mod}%
\ 3)$. A multigraph that satisfies these conditions is called $K_{3}%
$-\emph{divisible}. Kirkman showed that being $K_{3}$-divisible is also
sufficient for a complete graph to possess a triangle decomposition. A natural
question, therefore, concerns the density of non-complete triangle
decomposable graphs. Some work on this topic concerns a conjecture due to
Nash-Williams \cite{Nash}. A graph $G$ of order $n$ and minimum degree
$\delta(G)$ is $(1-\varepsilon)$-\emph{dense} if $\delta(G)\geq(1-\varepsilon
)(n-1)$. Nash-Williams conjectured that any sufficiently large $K_{3}%
$-divisible $\frac{3}{4}$-dense graph is $K_{3}$-decomposable. Keevash
\cite{Keev} obtained an asymptotic result, a special case of which applies to
this conjecture, with a value of $\varepsilon$ much smaller than $\frac{1}{4}$.

Holyer \cite{Holyer} showed that the problem of deciding whether a given
general graph is $K_{n}$-decomposable is NP-complete for $n\geq3$. Conditions
for different classes of planar graphs to be decomposable into paths of length
3 are presented in \cite{Hagg}. For decompositions of graphs into other graphs
$H$ of size $|E(H)|=3$, see e.g.~\cite{Rod, Fav, Lone}. On a somewhat
different note, planar graphs decomposable into a forest and a matching are
considered in several publications, including \cite{Borodin, Wang}, while it
is shown in \cite{Kim} that any planar graph is decomposable into three
forests, one of which has maximum degree at most four.

In contrast to the asymptotic results on $K_{n}$-decompositions of dense
graphs, we consider planar multigraphs and, in Section \ref{SecTriangle},
characterize those that are triangle decomposable. We begin with some
definitions and the statement of the characterization in Section
\ref{SecDefs}, followed by a number of lemmas in Section \ref{SecLemmas} and
the proof in Section \ref{SecProof}. In Section \ref{SecRational} we turn to
rational decompositions of planar multigraphs. We show in Section
\ref{SecWeights} that any rationally $K_{3}$-decomposable (simple) graph
admits such a decomposition using only weights $0,1$ or $\frac{1}{2}$, a
result which leads to a characterization of such graphs. We characterize
$K_{3}$-decomposable and rationally $K_{3}$-decomposable multigraphs that have
$K_{4}$ as underlying graph in Section \ref{SecK4}. We close with some ideas
for further work in Section~\ref{SecOpen}.

\section{Triangle Decompositions of Planar Multigraphs \label{SecTriangle}}

\subsection{Definitions and statement of main result}

\label{SecDefs}Since a multigraph is $K_{3}$-decomposable if and only if each
of its blocks is $K_{3}$-decomposable, we consider only 2-connected planar
multigraphs. In addition to being $K_{3}$-divisible, a $K_{3}$-decomposable
multigraph also needs to satisfy the condition that each of its edges is
contained in a triangle, a condition that holds trivially for (large enough)
complete graphs. A $K_{3}$-divisible multigraph that satisfies this third
necessary condition is called \emph{strongly} $K_{3}$-\emph{divisible}%
.\emph{\ }The planar graph $H$ obtained by joining the two vertices of
$K_{2,7}$ of degree seven shows that a strongly $K_{3}$-divisible graph need
not be $K_{3}$-decomposable: the removal of any triangle of $H$ results in a
triangle-free graph.%

\begin{figure}[ptb]%
\centering
\includegraphics[
height=2.399in,
width=2.2278in
]%
{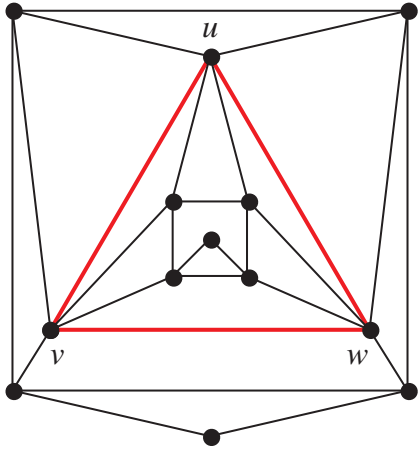}%
\caption{Triangle $uvw$ is a faceless triangle}%
\label{FigFaceless}%
\end{figure}

We denote a triangle with vertex set $\{u,v,w\}$ by $\tau=uvw$ if we are not
interested in the specific edges between its vertices. If specific edges are
important, we denote $\tau$ by $efg$, where $e=uv$, $f=vw$, and $g=wu$. A
triangle $\tau$ of a planar multigraph $G$ is called \emph{faced} if there
exists a plane embedding $\widetilde{G}$ of $G$ such that $\tau$ is a face of
$\widetilde{G}$; otherwise $\tau$ is called \emph{faceless}. The triangle
$uvw$ of the graph in Fig.~\ref{FigFaceless} is a faceless triangle; this can
be seen without much effort, but also follows from Lemma \ref{LemFaceless}
below. A \emph{separating triangle} $uvw$ of $G$ is one such that
$G-\{u,v,w\}$ is disconnected.

For vertices $u,v\in V(G)$, denote the number of edges joining $u$ and $v$ by
$\mu(u,v)$. A \emph{duplicate triangle} is a triangle $u_{1}u_{2}u_{3}$ such
that $\mu(u_{i},u_{j})\geq2$ for each $i\neq j$, and may be faced or faceless,
separating or non-separating. By \emph{deleting the edges of a duplicate
triangle} we mean that we delete exactly one edge between each pair of
vertices $u_{i}$ and $u_{j}$ of a duplicate triangle $u_{1}u_{2}u_{3} $.

A \emph{triangle depletion}, or simply a \emph{depletion}, of $G$ is any
spanning subgraph of $G$ obtained by sequentially deleting edges of (any
number of) faceless or duplicate triangles; note that $G$ is a depletion of itself.

The \emph{dual multigraph }$G^{\ast}$ of a plane multigraph $G$ is a plane
multigraph having a vertex for each face of $G$. The edges of $G^{\ast}$
correspond to the edges of $G$ as follows: if $e$ is an edge of $G$ that has a
face $F$ on one side and a face $F^{\prime}$ on the other side, then the
corresponding dual edge $e^{\ast}\in E(G^{\ast})$ is an edge joining the
vertices $f$ and $f^{\prime}$ of $G^{\ast}$ that correspond to the faces $F $
and $F^{\prime}$ of $G$. Note that under our assumption that $G$ is
2-connected, $G^{\ast}$ has no loops, and, using a careful geometric
description of the placement of vertices and edges in the dual, as in
\cite[Remark 7.1.8]{West}, we see that $(G^{\ast})^{\ast}\cong G$.

The statement of the main result of this section follows.

\begin{theorem}
\label{ThmTriangleDecomp}A planar multigraph $G$ is triangle decomposable if
and only if some depletion of $G$ has a plane embedding whose dual is a
bipartite multigraph in which all vertices of some partite set have degree three.
\end{theorem}

\subsection{Lemmas}

\label{SecLemmas}In our first result we present a characterization of faceless
triangles of planar multigraphs.\label{Fix w.r.t. multiple triangle deletion}

\begin{lemma}
\label{LemFaceless}A triangle $\tau=v_{1}v_{2}v_{3}$ of a planar multigraph
$G$ is faceless if and only if there exist two components $H_{1}$ and $H_{2}$
of $G-\{v_{1},v_{2},v_{3}\}$ such that each $v_{i}$ is adjacent, in $G$, to a
vertex in each $H_{j}$, $i=1,2,3,\ j=1,2$.
\end{lemma}

\noindent\textbf{Proof}.\hspace{0.1in}Let $\widetilde{G}$ be a plane embedding
of $G$ having $\tau$ as a face, but $G-\{v_{1},v_{2},v_{3}\}$ has components
$H_{j}$ as described. Let $G^{\prime}$ be the multigraph obtained by joining a
new vertex $v$ to each $v_{i}$. By inserting $v$ in the face $\tau$ of
$\widetilde{G}$, we get a plane embedding of $G^{\prime} $. However, by
contracting each $H_{i}$ to a single vertex we now obtain a $K_{3,3}$ minor of
$G^{\prime}$, a contradiction.

Conversely, suppose two such components $H_{j}$ do not exist. Let
$\widetilde{G}$ be a plane embedding of $G$ and suppose $\tau$ is not a face
of $\widetilde{G}$. Then $\widetilde{G}$ has vertices interior and exterior to
$\tau$. By assumption we may assume without loss of generality that each
component of $G-\{v_{1},v_{2},v_{3}\}$ interior to $\tau$ has vertices
adjacent, in $G$, to at most two vertices $v_{i},\ i=1,2,3$. Let $H$ be a
component of $G-\{v_{1},v_{2},v_{3}\}$ interior to $\tau$ such that no vertex
of $H$ is adjacent to (say) $v_{3}$. Let $F$ be the face of $\widetilde{G}$
exterior to $\tau$ that contains $v_{1}v_{2}$ on its boundary. By moving $H$
to $F$ we obtain an embedding of $G$ such that $H$ is exterior to $\tau$. By
repeating this procedure we eventually obtain an embedding $\widetilde{G}%
^{\prime}$ of $G$ such that $\tau$ is a face of $\widetilde{G}^{\prime}%
$.$~\blacksquare$

\bigskip

Evidently, then, a faceless triangle is a separating triangle.

\begin{lemma}
\label{Lem2_Conn}If a planar multigraph $G$ is $2$-connected, then so is any
depletion of $G$.
\end{lemma}

\noindent\textbf{Proof}.\hspace{0.1in}Suppose the statement of the lemma does
not hold, and let $G$ be a $2$-connected planar multigraph with the minimum
number of edges such that a depletion of $G$ is not $2$-connected. Then there
exists a faceless or duplicate triangle $\tau=uvw$ whose edges can be deleted
from $G$ to obtain a planar multigraph $G^{\prime}$ that is not $2$-connected.
This is impossible if $\tau$ is a duplicate triangle, hence $\tau$ is a
faceless triangle. Some vertex, say $v$, of $\tau$ is a cut-vertex of
$G^{\prime}$ but not of $G$.

Let $H$ be a component of $G-\{u,v,w\}$ whose existence is guaranteed by Lemma
\ref{LemFaceless}. Then both $u$ and $w$ are adjacent, in $G^{\prime}-v$, to
vertices of $H$. Therefore $u$ and $w$ belong to the same component, say $A$,
of $G^{\prime}-v$. Let $B$ be the union of all other components of $G^{\prime
}-v$. Then no vertex of $A$ is adjacent, in $G^{\prime}-v$, to a vertex of
$B$. Reinserting the edge $uw$ in $A$, we see that no vertex of $A+uw$ is
adjacent, in $G-v$, to a vertex of $B$; that is, $v$ is also a cut-vertex of
$G$, a contradiction.$~\blacksquare$

\bigskip

We also need the following result.

\begin{proposition}
\emph{\cite[Theorem 7.1.13]{West}} \label{PropDual}A plane multigraph is
Eulerian if and only if its dual is bipartite.
\end{proposition}

\subsection{Proof of Theorem \ref{ThmTriangleDecomp}}

\label{SecProof}We restate the characterization of triangle decomposable
planar multigraphs for convenience.

\bigskip

\noindent\textbf{Theorem }\ref{ThmTriangleDecomp}\hspace{0.1in}\emph{A planar
multigraph }$G$\emph{\ is triangle decomposable if and only if some depletion
}$G_{\Delta}$\emph{\ of }$G$\emph{\ has a plane embedding whose dual is a
bipartite multigraph in which all vertices of some partite set have degree
three.}

\bigskip

\noindent\textbf{Proof}.\hspace{0.1in}We may assume that $G$ is 2-connected.
Suppose $G$ is triangle decomposable. Then $G$ is strongly $K_{3}$-divisible.
Let $\mathcal{S}$ be the collection of triangles in some triangle
decomposition of $G$ and let $\mathcal{S}^{\prime}$ consist of all faceless
triangles, or triangles forming part of duplicate triangles, in $\mathcal{S}$.
Since the triangles in $\mathcal{S}^{\prime}$ are pairwise edge-disjoint,
deleting their edges results in a depletion $G_{\Delta}$ of $G$. Since
$\mathcal{S}$ is a triangle decomposition of $G$, $\mathcal{S-S}^{\prime}$ is
a triangle decomposition of $G_{\Delta}$, and every vertex of $G_{\Delta}$ is even.

Among all plane embeddings of $G_{\Delta}$, let $\widetilde{G}_{\Delta}$ be
one that maximizes the number of triangles in $\mathcal{S-S}^{\prime}$ that
are faces of the embedding. Suppose $\tau=uvw$ is a triangle in $\mathcal{S-S}%
^{\prime}$ that is not a face of $\widetilde{G}_{\Delta}$. Since $\tau$ is a
faced triangle, Lemma \ref{LemFaceless} implies that we may assume without
loss of generality that each component of $G-\{u,v,w\}$ interior to $\tau$ has
vertices adjacent, in $G$, to at most two of $u$, $v$ and $w$. Since $\tau$ is
not a duplicate triangle of $G_{\Delta}$, we may further assume that there is
at least one component of $G-\{u,v,w\}$ interior to $\tau$. Let $H$ be such a
component; say no vertex of $H$ is adjacent to $w$. Let $F$ and $F^{\prime}$
be the faces interior and exterior to $\tau$, respectively, containing the
edge $uv$ on their boundaries. Then neither $F$ nor $F^{\prime}$ is contained
in $\mathcal{S}$. By moving $H$ from $F$ to $F^{\prime}$ we obtain an
embedding of $G_{\Delta}$ such that $H$ is exterior to $\tau$. By repeating
this procedure we eventually obtain an embedding $\widetilde{G}_{\Delta
}^{\prime}$ of $G$ such that $\tau$ is a face of $\widetilde{G}_{\Delta
}^{\prime}$ and such that each triangle in $\mathcal{S-S}^{\prime}$ that is a
face of $\widetilde{G}_{\Delta}$ is also a face of $\widetilde{G}_{\Delta
}^{\prime} $. This contradicts the choice of $\widetilde{G}_{\Delta}$.

Hence all triangles in $\mathcal{S-S}^{\prime}$ are faces of $\widetilde{G}%
_{\Delta}$. By Lemma \ref{Lem2_Conn}, $G_{\Delta}$ is 2-connected. Thus each
edge of $\widetilde{G}_{\Delta}$ lies on two faces. Since $G_{\Delta}$ is
Eulerian, the dual $G_{\Delta}^{\ast}$ of $\widetilde{G}_{\Delta}$ is
bipartite (Proposition \ref{PropDual}). Let $(A,B)$ be a bipartition of
$G_{\Delta}^{\ast}$. Let $\tau,\tau^{\prime}$ be two triangles in
$\mathcal{S-S}^{\prime}$, let $t,t^{\prime}$ be the corresponding vertices of
$G_{\Delta}^{\ast}$ and assume without loss of generality that $t\in A$.
Consider any $t-t^{\prime}$ path $t=t_{0},t_{1},...,t_{k}=t^{\prime}$ in
$G_{\Delta}^{\ast}$ and say $t_{i}$ corresponds to a face $F_{i}$ of
$\widetilde{G}_{\Delta}$, $i=1,...,k$. Then $F_{1}$ is adjacent to $\tau$,
hence $F_{1}\notin\mathcal{S}$. Since $F_{2}$ is adjacent to $F_{1}$ and the
shared edge on the boundaries of $F_{1}$ and $F_{2}$ belongs to a triangle in
$\mathcal{S-S}^{\prime}$, $F_{2}\in\mathcal{S-S}^{\prime}$. Continuing this
argument we see that $F_{i}\in\mathcal{S-S}^{\prime}$ if and only if $i$ is
even. Since $F_{k}=\tau^{\prime}\in\mathcal{S-S}^{\prime}$, $k$ is even.
Therefore $t^{\prime}=t_{k}\in A$. We conclude that $A$ consists of all
vertices of $G_{\Delta}^{\ast}$ that correspond to triangles in $\mathcal{S-S}%
^{\prime}$, while all other vertices of $G_{\Delta}^{\ast}$ correspond to
faces of $\widetilde{G}_{\Delta}$ that are adjacent to triangles in
$\mathcal{S-S}^{\prime}$; hence these vertices belong to $B$. Therefore $\deg
v=3$ for all $v\in A$.

Conversely, suppose some depletion $G_{\Delta}$ of $G$ has a plane embedding
$\widetilde{G}_{\Delta}$ whose dual $G_{\Delta}^{\ast}$ possesses the stated
properties. By Proposition \ref{PropDual}, $G_{\Delta}$ is Eulerian. Let
$\mathcal{S}^{\prime}$ be the collection of edge disjoint triangles of $G$
whose deletion resulted in $G_{\Delta}$. Let $(A,B)$ be a bipartition of
$G_{\Delta}^{\ast}$ such that all vertices in $A$ have degree three and let
$\mathcal{S}$ be the faces of $\widetilde{G}_{\Delta}$ corresponding to the
vertices in $A$. Since $A$ is an independent set of vertices that cover all
edges of $G_{\Delta}^{\ast}$ (since $G_{\Delta}^{\ast}$ is a multigraph, it
has no loops), $\mathcal{S}$ consists of mutually edge-disjoint triangles
covering all edges of $G_{\Delta}$. Therefore $\mathcal{S}$ is a triangle
decomposition of $G_{\Delta}$ and $\mathcal{S}\cup\mathcal{S}^{\prime}$ is a
triangle decomposition of $G$.~$\blacksquare$

\bigskip

\begin{figure}[ptb]%
\centering
\includegraphics[
height=2.3912in,
width=4.9286in
]%
{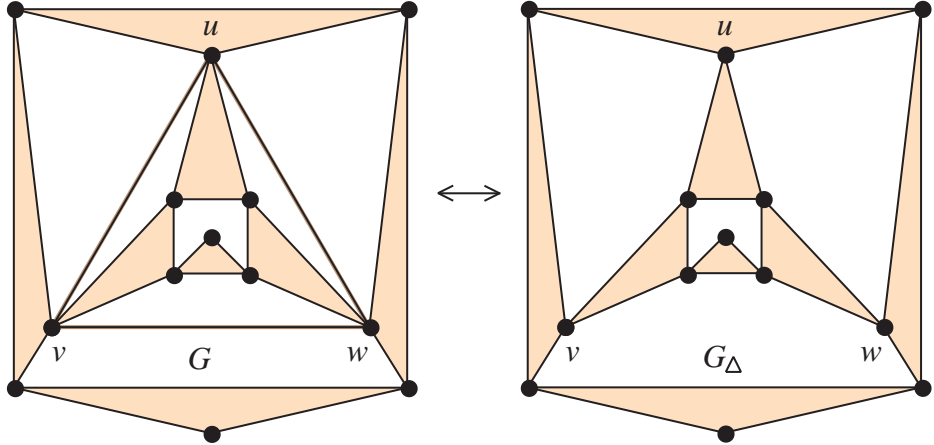}%
\caption{Triangle decompositions of $G$ and $G_{\Delta}$}%
\label{FigDecomp}%
\end{figure}

Triangle decompositions of a graph $G$ and its depletion $G_{\Delta}$ are
illustrated in Fig.~\ref{FigDecomp}. Since $G$ itself is Eulerian, the dual of
any embedding of $G$ is bipartite. However, no embedding of $G$ has a dual in
which all vertices of one partite set of its bipartition have degree three:
the edge $vw$ always lies on two nontriangular faces, and the corresponding
vertices (of degree at least four) of the dual are in different partite sets.
A $K_{3}$-decomposition of $G$ is obtained by first deleting $uvw$,
partitioning the faces into two sets so that one set contains only triangles,
which form part of the decomposition, and reinserting $uvw$ to complete the decomposition.

Theorem \ref{ThmTriangleDecomp} implies that the necessary conditions for a
multigraph to be triangle decomposable are also sufficient for maximal planar
graphs, which trivially satisfy two of the conditions (of being strongly
$K_{3}$-divisible) provided they have order at least three.%

\begin{corollary}
A maximal planar graph is triangle decomposable if and only if it is Eulerian.
\end{corollary}

\noindent\textbf{Proof}.\hspace{0.1in}Any plane embedding of a maximal planar
graph $G$ of order at least three is a triangulation of the plane. Its dual is
cubic, and bipartite because $G$ is Eulerian, and either partite set
corresponds to a triangle decomposition of $G$.~$\blacksquare$

\begin{corollary}
Any Eulerian multigraph whose edges can be partitioned into sets that induce
maximal planar subgraphs is triangle decomposable.
\end{corollary}

\section{Rational Triangle Decompositions}

\label{SecRational}The main purpose of this section is to characterize
rationally triangle decomposable planar graphs, which we do in Corollary
\ref{CorRat}, after first showing, in Theorem \ref{ThmRational}, that each
such graph admits a rational triangle decomposition using only weights $0,1$
and $\frac{1}{2}$. In Section \ref{SecK4} we characterize $K_{3}$-decomposable
and rationally $K_{3}$-decomposable planar multigraphs that have $K_{4}$ as
underlying graph in terms of the multiplicities of their edges.

Dense graphs that admit rational $K_{3}$-decompositions were studied in
\cite{Kseniya}. The only condition among the three for a multigraph $G$ to be
$K_{3}$-decomposable that remains necessary for $G$ to be rationally $K_{3}%
$-decomposable is the condition that each edge of $G$ be contained in a
triangle. Clearly, maximal planar graphs of order at least three are
rationally $K_{3}$-decomposable: assign a weight of $\frac{1}{2}$ to each face
triangle in a plane embedding of the graph. In fact, each multigraph whose
edges can be partitioned into sets that induce maximal planar subgraphs is
rationally $K_{3}$-decomposable.

\subsection{Rationally triangle decomposable planar graphs}

\label{SecWeights}Suppose $G$ is a rationally $K_{3}$-decomposable multigraph
and consider such a decomposition of $G$. For a triangle $\tau$ of $G$, we
denote the weight of $\tau$ by $w(\tau)$, and for any edge $e$ of $G$, we
denote the sum of the weight of the triangles that contain $e$ by $w(e)$;
since $G$ is rationally $K_{3}$-decomposable, $w(e)=1$ for each edge $e$.

While it is easy to find planar multigraphs that possess rational triangle
decompositions with weights $\frac{p}{q}$ and $\frac{q-p}{q}$ for arbitrary
integers $q\geq1$ and $0\leq p\leq q$, for example Eulerian maximal planar
graphs, all examples of rationally $K_{3}$-decomposable multigraphs we know of
also admit decompositions using only weights $0,\ 1$ and $\frac{1}{2}$. We
show that this is true for all rationally $K_{3}$-decomposable (simple) planar graphs.

For a triangle $\tau=xyz$ of a plane graph $G$, let $I_{\tau}$ denote the
subgraph of $G$ induced by $\{x,y,z\}$ and all vertices interior to $\tau$. We
call $I_{\tau}$ the \emph{interior graph} of $\tau$. A separating triangle of
$G$ is an \emph{innermost }(or an \emph{outermost}) \emph{separating triangle}
if its interior (or its exterior) contains no separating triangles. Similarly,
a separating triangle containing an edge $e$ is an \emph{outermost separating
triangle containing} $e$ if no separating triangle in its exterior
contains$~e$.

\begin{theorem}
\label{ThmRational}If $G$ is a rationally $K_{3}$-decomposable planar graph,
then $G$ has a $K_{3}$-decompo-sition using only weights $0$, $1$, and
$\frac{1}{2}$.
\end{theorem}

\noindent\textbf{Proof}.\hspace{0.1in}Suppose there exists a planar graph that
is rationally $K_{3}$-decomposable but does not have a decomposition using
only weights $0$, $1$, and $\frac{1}{2}$. Let $H$ be such a graph with the
minimum number of edges. We establish the following properties of $H$:

\begin{enumerate}
\item $H$\emph{\ is not maximal planar:\hspace{0.1in}}A maximal planar graph
has a $K_{3}$-decomposition where each face receives weight $\frac{1}{2}$.

\item \emph{Every edge of }$H$\emph{\ is in at least two triangles:}%
\textbf{\hspace{0.1in}}Suppose $e\in E(H)$ is in only one triangle $\tau$.
Then in any rational $K_{3}$-decomposition of $H$, $\tau$ receives weight 1.
Hence $H-\tau$ is rationally $K_{3}$-decomposable, and since $H-\tau$ has
fewer edges than $H$, it has a rational $K_{3}$-decomposition using only
weights $0 $, $1$, and $\frac{1}{2}$. But then $H$ has a decomposition using
only weights $0$, $1$, and $\frac{1}{2}$, which is a contradiction.

\item \emph{In any embedding of }$H$\emph{, every edge incident with a
nontriangular face belongs to a separating triangle:\hspace{0.1in}}Let $e$ be
an edge incident with a nontriangular face. As $e$ is in at least two
triangles, and is incident with exactly two faces, one such triangle $\tau$ is
not a face. Thus there are vertices interior and exterior to $\tau$, which is
therefore a separating triangle.

\item \emph{In any embedding of }$H$\emph{, there exists a separating triangle
}$\tau$,\emph{\ incident with an edge }$e$\emph{\ of a nontriangular face,
whose exterior contains no triangles containing }$e$\emph{\ and whose interior
graph }$I_{\tau}$\emph{\ is maximal planar:\hspace{0.1in}}Let $e_{1}$ be an
edge of a nontriangular face and let $\tau_{1}$ be the outermost separating
triangle containing $e$. If $I_{\tau_{1}}$ is maximal planar, we are done.
Otherwise, $I_{\tau_{1}}$ contains a nontriangular face; choose an edge
$e_{2}$ of this face not also contained in $\tau_{1}$ and its outermost
separating triangle $\tau_{2}$. Note that $\tau_{2}$ lies interior to
$\tau_{1}$. As $H$ is finite, this process terminates.
\end{enumerate}

Now, assume that $H$ is embedded in the plane and that $\tau$ is a separating
triangle incident with an edge $e$ of a nontriangular face $f$, whose exterior
contains no triangles containing $e$ and whose interior graph $I_{\tau}$ is
maximal planar. We next prove the following claim regarding an innermost
separating triangle of $H$.

\bigskip

\noindent\textbf{Claim \ref{ThmRational}.1}\hspace{0.1in}\textit{Let }%
$T$\textit{\ be an innermost separating triangle of }$H$\textit{. Then }%
$I_{T}$\textit{\ is maximal planar and any rational }$K_{3}$%
\textit{-decomposition of }$H$\textit{\ gives the same weight to the faces of
}$I_{T}$\textit{\ adjacent to }$T$\textit{.}\bigskip

\noindent\textbf{Proof of Claim \ref{ThmRational}.1}.\hspace{0.1in}Suppose
$I_{T}$ is not maximal planar. Then it contains a nontriangular face. But
every edge of this face that does not belong to $T$ is in a separating
triangle interior to $T$, which contradicts the choice of $T$. Hence $I_{T}$
is maximal planar.

Let $I_{T}^{\ast}$ be the dual of $I_{T}$ and let $D=I_{T}^{\ast}-T$. First
suppose $D$ is bipartite with bipartition $(X,Y)$. We show that the vertices
$x,\ y$ and $z$ of $D$ corresponding to the three faces of $I_{T}$ adjacent to
$T$ are in the same partite set. Otherwise, assume without loss of generality
that $x,y\in X$ and $z\in Y$. Since $I_{T}$ is maximal planar, $x,\ y$ and $z$
have degree 2 and every other vertex of $D$ has degree 3. But then the number
of edges incident with a vertex in $X$ is congruent to $1\ (\operatorname{mod}%
\ 3)$ and the number of edges incident with a vertex in $Y$ is congruent to
$2\ (\operatorname{mod}\ 3)$, which is impossible. Hence $x,\ y$ and $z$
belong to the same partite set.

Now, since every edge in $D$ corresponds to an edge of $I_{T}$ that lies on
exactly two triangle faces, and $I_{T}$ has no separating triangles, every
face in the same partite set receives the same weight, and the weights of the
two sets sum to 1. Hence, the faces of $I_{T}$ adjacent to $T$ receive the
same weight.

Now suppose $D$ is not bipartite. Then $D$ contains an odd cycle $f_{1}%
f_{2}f_{3}\cdots f_{k}f_{1}$. Suppose $w(f_{1})=x$. Then as every edge is
incident with exactly two triangles, $w(f_{2})=1-x$, $w(f_{3})=x$, \ldots,
$w(f_{k})=x$, and $w(f_{1})=1-x=x$. Hence $x=\frac{1}{2}$. Filling in the
remaining weights from this cycle, every face receives weight $\frac{1}{2}$.
Hence, the faces of $I_{T}$ adjacent to $T$ receive the same weight.~$\square$

\bigskip

Continuing with the proof of Theorem \ref{ThmRational}, consider a rational
$K_{3}$-decomposition of $H$. Suppose $I_{\tau}$ contains a separating
triangle other than $\tau$. Then choose an innermost separating triangle
$\tau^{\prime}$ of $I_{\tau}$ and let $H^{\prime}$ be the graph obtained by
deleting the interior of $\tau^{\prime}$. By Claim \ref{ThmRational}.1, the
interior faces of $I_{\tau^{\prime}}$ adjacent to $\tau^{\prime}$ receive the
same weight, say $x$. Then the rational $K_{3}$-decomposition of $H$ induces a
rational $K_{3}$-decomposition of $H^{\prime}$ in which $w_{H^{\prime}}%
(\tau^{\prime})=w_{H}(\tau^{\prime})+x$. We continue this process until $\tau$
has no separating triangles in its interior. Finally, apply this process to
$\tau$ itself, obtaining the graph $H^{\dagger}$. Now $e$ is contained in only
one triangle in $H^{\dagger}$, namely $\tau$, so $w_{H^{\dagger}}(\tau)=1$.
Then $H^{\dagger}-\tau$ has a rational $K_{3}$-decomposition and since
$H^{\dagger}-\tau$ has fewer edges than $H$, it has a rational $K_{3}%
$-decomposition using only weights $0$, $1$, and $\frac{1}{2}$. As a result,
we obtain a rational $K_{3}$-decomposition of $H$ using only weights $0$, $1$,
and $\frac{1}{2}$ by extending the decomposition of $H^{\dagger}-\tau$ and
giving each face of $I_{\tau}$ (including $\tau$) weight $\frac{1}{2}$. This
decomposition contradicts the assumption that $H$ does not have a
decomposition using only weights $0$, $1$, and $\frac{1}{2}$, completing the
proof.~$\blacksquare$

\bigskip

Let $^{2}G$ denote the multigraph obtained from a simple graph $G$ by
replacing each edge by a pair of parallel edges. For an edge $e$ of $G$,
denote the corresponding pair of edges of $^{2}G$ by $e_{1}$ and $e_{2}$. If
$\tau_{1}$ and $\tau_{2}$ are edge disjoint triangles of $^{2}G$ with the same
vertex set, denote the corresponding triangle of $G$ by $\tau$. For $u,v\in
V(^{2}G)$, denote the set of edges joining $u$ and $v$ by $E(u,v)$. The
characterization of rationally triangle decomposable planar graphs follows.

\begin{corollary}
\label{CorRat}A simple planar graph $G$ is rationally $K_{3}$-decomposable if
and only if $^{2}G$ is $K_{3}$-decomposable.
\end{corollary}

\noindent\textbf{Proof.\hspace{0.1in}}Suppose $G$ is rationally $K_{3}%
$-decomposable. By Theorem \ref{ThmRational}, $G$ has a rational $K_{3}%
$-decomposition using only weights $0$, $1$, and $\frac{1}{2}$. Let
$\mathcal{T}_{\frac{1}{2}}$ and $\mathcal{T}_{1}$ denote the sets of triangles
of $G$ with weights $\frac{1}{2}$ and $1$, respectively. For any triangle
$efg\in\mathcal{T}_{1}$, partition the edges $e_{i},f_{i},g_{i},\ i=1,2$, of
$^{2}G$ arbitrarily into two triangles $\tau_{1}$ and $\tau_{2}$, and let
$w(\tau_{1})=w(\tau_{2})=1$. Any edge of $G$ that belongs to a triangle in
$\mathcal{T}_{\frac{1}{2}}$ belongs to exactly two triangles in $\mathcal{T}%
_{\frac{1}{2}}$ and to no triangles in $\mathcal{T}_{1}$. Therefore, for the
set of edges of $G$ that belong to triangles in $\mathcal{T}_{\frac{1}{2}}$,
the corresponding set of edge pairs of $^{2}G$ can be partitioned into edge
disjoint triangles, each being allocated weight $1$, to give a $K_{3}%
$-decomposition of $^{2}G$.

Conversely, assume $^{2}G$ is $K_{3}$-decomposable. For vertices $x,y,z$ of
$^{2}G$ and edges $e_{1},e_{2}\in E(x,y),\ f_{1},f_{2}\in E(y,z)$ and
$g_{1},g_{2}\in E(x,z)$, if $e_{i},f_{i},g_{i},\ i=1,2$, can be partitioned
into triangles $\tau_{1}$ and $\tau_{2}$ such that $w(\tau_{1})=w(\tau_{2}%
)=1$, let $w(\tau)=1$, and if $e_{i},f_{i},g_{i}$ can be partitioned into
triangles $\tau_{1}$ and $\tau_{2}$ such that (say) $w(\tau_{1})=0$ and
$w(\tau_{2})=1$, let $w(\tau)=\frac{1}{2}$. Since each edge that belongs to
$\tau_{1}$ also belongs to another triangle of $^{2}G$ with weight $1$, this
gives a rational $K_{3}$-decomposition of $G$.~$\blacksquare$

\begin{corollary}
\label{Cor0mod3}If $G$ is a rationally $K_{3}$-decomposable planar graph, then
$|E(G)|\equiv0\ (\operatorname{mod}\ 3)$.
\end{corollary}

\noindent\textbf{Proof.\hspace{0.1in}}By Corollary \ref{CorRat}, $^{2}G$ has a
$K_{3}$-decomposition. Therefore $|E(^{2}G)|\equiv0\ (\operatorname{mod}\ 3)$.
Since $|E(^{2}G)|=2|E(G)|$, we also have $|E(G)|\equiv0\ (\operatorname{mod}%
\ 3)$.~$\blacksquare$

\subsection{Multigraphs with $K_{4}$ as underlying graph}

\label{SecK4}One reason why the proof of Theorem \ref{ThmRational} fails for
multigraphs is that multiple edges that do not lie on triangular faces are not
necessarily contained in separating triangles. Hence statement (4) in the
proof does not necessarily hold; certainly, if its underlying graph is
complete, a multigraph contains no separating triangles at all.

It is easy to see that a multigraph $G$ with $K_{3}$ as underlying graph is
rationally $K_{3}$-decomposable if and only if all edges have the same
multiplicity, say $k$; in this case, $|E(G)|=3k$ and $G$ can be decomposed
into $k$ edge-disjoint triangles. In the remainder of this section we
characterize $K_{3}$-decomposable and rationally $K_{3}$-decomposable
multigraphs that have $K_{4}$ as underlying graph.

Denote the set of all multigraphs that have $K_{4}$ as underlying graph by
$\mathcal{K}_{4}$. For any $G\in\mathcal{K}_{4}$ and distinct edges $e$ and
$f$, let $w(e,f)$ be the sum\ of the weight of the triangles that contain both
$e$ and $f$, and for any vertices $u,v$ of $G$, let $w(uv,e)$ be the sum of
the weight of the triangles that contain $e$ and some edge joining $u$ and
$v$. Also, for $u,v\in V(G)$, denote the set of edges joining $u$ and $v$ by
$E(u,v)$.

\begin{figure}[ptb]%
\centering
\includegraphics[
height=2.719in,
width=2.8746in
]%
{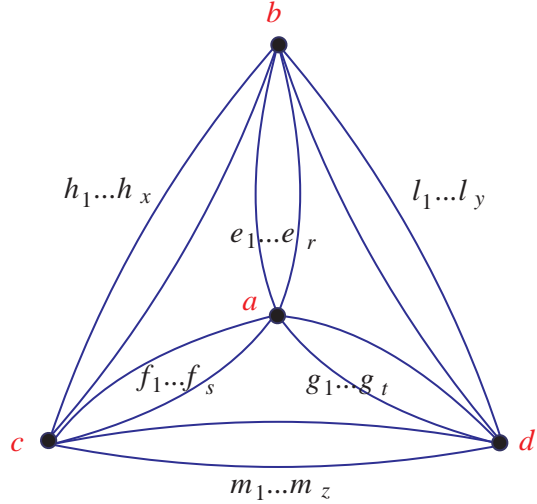}%
\caption{Labels of the vertices and edges of the multigraph $G$ with $K_{4}$
as underlying graph}%
\label{FigK_4}%
\end{figure}

Say $V(G)=\{a,b,c,d\}$. The following notation will be used throughout this
subsection (see Fig.~\ref{FigK_4}). Let $\mu(a,b)=r,\ \mu(a,c)=s,\ \mu
(a,d)=t,\ \mu(b,c)=x,\ \mu(b,d)=y$ and $\mu(c,d)=z$, and let
\begin{equation}%
\begin{tabular}
[c]{l}%
$E(a,b)=\{e_{1},...,e_{r}\}$\\
$E(a,c)=\{f_{1},...,f_{s}\}$\\
$E(a,d)=\{g_{1},...,g_{t}\}$\\
$E(b,c)=\{h_{1},...,h_{x}\}$\\
$E(b,d)=\{\ell_{1},...,\ell_{y}\}$\\
$E(c,d)=\{m_{1},...,m_{z}\}.$%
\end{tabular}
\label{eqNotation}%
\end{equation}
%

\begin{theorem}
\label{Thm_K4}Let $G\in\mathcal{K}_{4}$, let $u\in V(G)$ and let
$V(G)-\{u\}=\{v_{1},v_{2},v_{3}\}$. Then $G$ is $K_{3}$-decomposable if and
only if

\begin{enumerate}
\item[$(i)$] there exists an integer $n$ such that $0\leq n\leq\min\{\mu
(v_{i},v_{j}):i,j\in\{1,2,3\},i\neq j\}$ and $\mu(u,v_{i})=\mu(v_{i}%
,v_{j})+\mu(v_{i},v_{k})-2n$ for each $i\in\{1,2,3\}$, each $j\in
\{1,2,3\}-\{i\}$ and $k\in\{1,2,3\}-\{i,j\}$,
\end{enumerate}

and rationally $K_{3}$-decomposable if and only if

\begin{enumerate}
\item[$(ii)$] there exists an integer $n^{\prime}$ such that $0\leq
\frac{n^{\prime}}{2}\leq\min\{\mu(v_{i},v_{j}):i,j\in\{1,2,3\},i\neq j\}$ and
$\mu(u,v_{i})=\mu(v_{i},v_{j})+\mu(v_{i},v_{k})-n^{\prime}$ for each
$i\in\{1,2,3\}$, each $j\in\{1,2,3\}-\{i\}$ and $k\in\{1,2,3\}-\{i,j\}$.
\end{enumerate}

Moreover, if $G$ is rationally $K_{3}$-decomposable, it has a decomposition
using only weights $0$, $1$ and $\frac{1}{2}$.
\end{theorem}

\noindent\textbf{Proof.}\hspace{0.1in}To simplify notation, let
$V(G)=\{a,b,c,d\}$ and assume without loss of generality that $(i)$ holds with
$u=a$. With notation as in (\ref{eqNotation}), if $n>0$, let
\begin{equation}
G^{\prime}=G-\bigcup_{i=0}^{n-1}\{h_{x-i}\ell_{y-i}m_{z-i}\}\label{eqG'1}%
\end{equation}
and let $x^{\prime}=x-n,\ y^{\prime}=y-n$ and $z^{\prime}=z-n$. Then in
$G^{\prime}$,
\[%
\begin{tabular}
[c]{l}%
$E^{\prime}(a,b)=\{e_{1},...,e_{r}\}$\\
$E^{\prime}(a,c)=\{f_{1},...,f_{s}\}$\\
$E^{\prime}(a,d)=\{g_{1},...,g_{t}\}$\\
$E^{\prime}(b,c)=\{h_{1},...,h_{x^{\prime}}\}$\\
$E^{\prime}(b,d)=\{\ell_{1},...,\ell_{y^{\prime}}\}$\\
$E^{\prime}(c,d)=\{m_{1},...,m_{z^{\prime}}\}$%
\end{tabular}
\]
and $(i)$ holds for $G^{\prime}$ and $a$ with $n=0$. Now%
\begin{equation}
E(G^{\prime})=\left(  \bigcup_{i=1}^{x^{\prime}}\{e_{i}f_{i}h_{i}\}\right)
\cup\left(  \bigcup_{i=1}^{y^{\prime}}\{e_{x^{\prime}+i}g_{i}\ell
_{i}\}\right)  \cup\left(  \bigcup_{i=1}^{z^{\prime}}\{f_{x^{\prime}%
+i}g_{y^{\prime}+i}m_{i}\}\right)  .\label{eqG'2}%
\end{equation}
Since each set of three edges in (\ref{eqG'1}) and (\ref{eqG'2}) induces a
triangle, $G$ is $K_{3}$-decomposable.

Conversely, suppose $G$ is $K_{3}$-decomposable into $\alpha$ triangles
induced by $\{b,c,d\}$, $\beta$ triangles induced by $\{a,c,d\}$, $\psi$
triangles induced by $\{a,b,d\}$ and $\delta$ triangles induced by
$\{a,b,c\}$. Then%
\begin{align*}
\mu(a,b)  & =\psi+\delta\\
\mu(a,c)  & =\beta+\delta\\
\mu(a,d)  & =\beta+\psi\\
\mu(b,c)  & =\alpha+\delta,\text{ hence }\delta=\mu(b,c)-\alpha\\
\mu(b,d)  & =\alpha+\psi,\text{ hence }\psi=\mu(b,d)-\alpha\\
\mu(c,d)  & =\alpha+\beta,\text{ hence }\beta=\mu(c,d)-\alpha
\end{align*}
and thus%
\begin{align*}
\mu(a,b)  & =\mu(b,c)+\mu(b,d)-2\alpha\\
\mu(a,c)  & =\mu(b,c)+\mu(c,d)-2\alpha\\
\mu(a,d)  & =\mu(b,d)+\mu(c,d)-2\alpha.
\end{align*}
Since $\beta,\psi,\delta\geq0$, $\alpha\leq\min\{\mu(b,c),\mu(b,d),\mu
(c,d)\}$. Therefore $(i)$ holds for $u=a$ and $n=\alpha$. Similarly, $(i)$
holds for $b,c$ and $d$ with $n=\beta,\psi$ and $\delta$, respectively.

\bigskip

Suppose $(ii)$ holds with $u=a$. If $n^{\prime}$ is even, let $n^{\prime}=2n$.
Then $(i)$ holds and $G$ is $K_{3}$-decomposable. Hence assume $n^{\prime}$ is
odd. Say $n^{\prime}=2n+1$ and let
\[
G^{\prime}=G-\{e_{r},f_{s},g_{t},h_{x},\ell_{y},m_{z}\}.
\]
Since $n+1\leq\min\{x,y,z\}$, $n\leq\min\{x-1,y-1,z-1\}$. The equations
$r=x+y-2n-1,\ s=x+z-2n-1$ and $t=y+z-2n-1$ in $G$ imply the equations
$r-1=(x-1)+(y-1)-2n,\ s-1=(x-1)+(z-1)-2n$ and $t-1=(y-1)+(z-1)-2n$ in
$G^{\prime}$. Hence $(i)$ holds for $G^{\prime}$ with $u=a$, and $G^{\prime}$
is $K_{3}$-decomposable.

Since $\{e_{r},f_{s},g_{t},h_{x},\ell_{y},m_{z}\}$ induces a $K_{4}$, which is
rationally $K_{3}$-decomposable into four triangles, each of weight $\frac
{1}{2}$, $G$ is rationally $K_{3}$-decomposable using only weights of $0 $,
$1$ and $\frac{1}{2}$.

Conversely, say $G$ is rationally $K_{3}$-decomposable and consider such a
decomposition of $G$. For each edge $e_{j}\in E(a,b)$, any triangle that
contains $e_{j}$ also contains one edge in $E(b,c)\cup E(b,d)$. Since
$w(e_{j})=1$,%
\[
\sum_{i=1}^{x}w(e_{j},h_{i})+\sum_{i=1}^{y}w(e_{j},\ell_{i})=1,
\]
hence%
\begin{equation}
\sum_{i=1}^{x}w(ab,h_{i})+\sum_{i=1}^{y}w(ab,\ell_{i})=r.\label{eq_ab1}%
\end{equation}
Similarly,%
\[
\sum_{i=1}^{x}w(ac,h_{i})+\sum_{i=1}^{z}w(ac,m_{i})=s
\]
and%
\[
\sum_{i=1}^{y}w(ad,\ell_{i})+\sum_{i=1}^{z}w(ad,m_{i})=t.
\]

Let $\mathcal{T}$ be the set of all triangles that do not contain any edges
incident with $a$, that is, triangles of the form $h_{i}\ell_{j}%
m_{k},\ i=1,...,x,\ j=1,...,y,\ k=1,...,z$, and let $\omega$ be the total
weight of the triangles in $\mathcal{T}$. Then $\omega\leq\min\{x,y,z\}$. For
any edge $h_{i}$, any triangle that contains $h_{i}$ but no edge in $E(a,b)$
belongs to $\mathcal{T}$. Hence $\sum_{i=1}^{x}w(ab,h_{i})+\omega=x$.
Similarly, $\sum_{i=1}^{y}w(ab,\ell_{i})+\omega=y$. Substitution in
(\ref{eq_ab1}) gives $r=x+y-2\omega$. Similarly, $s=x+z-2\omega$ and
$t=y+z-2\omega$. Since $r,x,y$ are integers, $2\omega$ is an integer, say
$2\omega=n^{\prime}$. Then $(ii)$ holds for $a$. As before, similar arguments
show that $(ii)$ also holds for $b,c$ and $d$.

As shown above, if $(ii)$ holds, then $G$ is rationally $K_{3}$-decomposable
using only weights of $0$, $1$ and $\frac{1}{2}$. This proves the last part of
the theorem.~$\blacksquare$

\bigskip

\label{here}By taking $\mu(v_{1},v_{2})=0$ in Theorem \ref{Thm_K4}$(ii)$, we
get the following corollary.

\begin{corollary}
\label{CorK4-e}Let $G$ be a multigraph whose underlying graph is $K_{4}-e$.
Say $V(G)=\{u,v,v_{1},v_{2}\}$, where $u$ and $v$ correspond to the vertices
of $K_{4}-e$ of degree three. The following conditions are equivalent:

\begin{enumerate}
\item $G$ is rationally $K_{3}$-decomposable.

\item $G$ is $K_{3}$-decomposable.

\item $\mu(u,v)=\mu(v,v_{1})+\mu(v,v_{2}),\ \mu(u,v_{1})=\mu(v,v_{1})$ and
$\mu(u,v_{2})=\mu(v,v_{2}).$
\end{enumerate}
\end{corollary}

The final corollary now follows similar to Corollary \ref{Cor0mod3}.

\begin{corollary}
If $G$ is a rationally $K_{3}$-decomposable multigraph whose underlying graph
is $K_{3},\ K_{4}$ or $K_{4}-e$, then $|E(G)|\equiv0\ (\operatorname{mod}\ 3)$.
\end{corollary}

\section{Open Questions}

\label{SecOpen}

\begin{enumerate}
\item Does Theorem \ref{ThmRational} hold for rationally $K_{3}$-decomposable
planar multigraphs?

\item Can we characterize rationally $K_{3}$-decomposable planar multigraphs?

\item Can we characterize rationally $K_{3}$-decomposable outerplanar graphs
or multigraphs?

\item What can we say about graphs embeddable on other surfaces?
\end{enumerate}

\end{document}